\documentclass{elsart}
\usepackage{amsmath}
\usepackage{amssymb}

\newcommand\im{\mathrm{im}}
\newcommand\tr{\mathrm{tr}}
\newcommand\ef{{\mathcal F}}
\newcommand\te{{\mathcal T}}
\newcommand\ce{{\mathcal C}}
\newcommand\mem{{\mathcal M}}
\newcommand\vng{{\mathcal N}G}
\newcommand\B{{\mathcal B}}
\newcommand\T{\mathrm{{\bf T}}}
\newcommand\bigP{\mathrm{{\bf P}}}
\newcommand\smallt{\mathrm{{\bf t}}}
\newcommand\p{\mathrm{{\bf p}}}
\newcommand\unb{\mathrm{{\bf u}}}
\newcommand\bnd{\mathrm{{\bf b}}}
\newcommand\cl{\mathrm{cl}}
\newcommand\homo{\mathrm{Hom}}
\newcommand\ann{\mathrm{ann}}

\begin{document}
\begin{frontmatter}

\title{Dimension and Torsion Theories for a Class of Baer *-Rings}

\author{Lia Va\v s}

\address{Department of Mathematics, Physics and Computer Science,
University of the Sciences in Philadelphia, 600 S. 43rd St.,
Philadelphia, PA 19104}

\ead{l.vas@usip.edu}

\begin{abstract}
Many known results on finite von Neumann algebras are generalized,
by purely algebraic proofs, to a certain class $\ce$ of finite
Baer *-rings. The results in this paper can also be viewed as a
study of the properties of Baer *-rings in the class $\ce$.

First, we show that a finitely generated module over a ring from
the class $\ce$ splits as a direct sum of a finitely generated
projective module and a certain torsion module. Then, we define
the dimension of any module over a ring from $\ce$ and prove that
this dimension has all the nice properties of the dimension
studied in \cite{Lu2} for finite von Neumann algebras. This
dimension defines a torsion theory that we prove to be equal to
the Goldie and Lambek torsion theories. Moreover, every finitely
generated module splits in this torsion theory.

If $R$ is a ring in $\ce,$ we can embed it in a canonical way into
a regular ring $Q$ also in $\ce.$ We show that $K_0(R)$ is
isomorphic to $K_0(Q)$ by producing an explicit isomorphism and
its inverse of monoids Proj$ (P)\rightarrow$ Proj$ (Q)$ that
extends to the isomorphism of $K_0(R)$ and $K_0(Q)$.
\end{abstract}


\end{frontmatter}

\section{Introduction}

This paper is motivated by the remark of Sterling K. Berberian
from the introduction to his book \cite{Be2} on Baer *-rings: "The
subject of Baer *-rings has its roots in von Neumann theory of
'rings of operators' (now called von Neumann algebras) ... Von
Neumann algebras are blessed with an excess of structure --
algebraic, geometric, topological -- so much, that one can easily
obscure, through proof by overkill, what makes a particular
theorem work." Relying just on algebra, we generalize some results
from \cite{Lu2}, \cite{Lia1} and \cite{Tele_Teza} to a certain
class of finite Baer-* rings that contains the class of all finite
von Neumann algebras. The proofs in \cite{Lu2} rely on some of the
geometric or topological structure of finite von Neumann algebras.
All proofs in this paper rely strictly on algebraic properties. We
follow Berberian's idea: "if all the functional analysis is
stripped away ... what remains should (be) completely accessible
through algebraic avenues".

We impose some restrictions onto the Baer *-ring $R$ that are
sufficient for defining the dimension function, the regular ring
of $R$, and ensuring that all the matrix algebras over $R$ are
sufficiently nice (we need the lattice of projections of all the
matrix algebras to be complete). In Section \ref{SectionClassCe},
we list the axioms imposed onto the Baer *-ring. These axioms are
the same ones Berberian uses in \cite{Be2} in order to ensure that
the matrix rings over $R$ are Baer *-rings. Baer *-rings that
satisfy those axioms form a class that we shall denote by $\ce$.
Every finite $AW^*$-algebra (so a finite von Neumann algebra in
particular) is in $\ce.$

As in \cite{Lia1}, torsion theories are used to study the modules
over the rings of interest. In Section \ref{SectionTT}, we recall
the definition of an arbitrary torsion theory and some related
notions. Then we list the examples of torsion theories we shall
use in the rest of the paper (Lambek, Goldie, classical, etc.).

In Section \ref{SectionTTforCe}, we prove the main results. First,
we show that a finitely generated module over a ring from class
$\ce$ splits as a direct sum of a finitely generated projective
module and a certain torsion module (Theorem
\ref{SplittingOfBnd}). This generalizes an analogous result for
finite von Neumann algebras proven in \cite{Lu2}.

Secondly, if $R$ is a ring from the class $\ce,$ we prove (Theorem
\ref{DimensionInR}) that a dimension of {\em any} $R$-module can
be defined so that it has all the nice properties of the dimension
defined in \cite{Lu2} (i.e. we prove that Theorem 0.6 from
\cite{Lu2} holds for the class $\ce$). This dimension defines a
torsion theory that proves to be equal to the Goldie and Lambek
torsion theories and every finitely generated module splits in
this torsion theory (Theorem \ref{dim=G=L}). In Theorem
\ref{TTEqual}, we demonstrate how the torsion theories reflect the
ring-theoretic properties of $R.$

If $R$ is a finite von Neumann algebra, our construction gives us
precisely the central-valued dimension considered in \cite{Lu1}
for finitely presented modules. Moreover, Theorem
\ref{DimensionInR}, guarantees that we can extend the definition
to any $R$-module. Thus, in the case of finite von Neumann
algebras, we can define the real-value dimension (as in \cite{Lu2}
or \cite{Lu_book}) and the central-valued dimension. In Section
\ref{SectionCorollaries}, we show that both dimensions define the
same torsion theory (Corollary \ref{CorollaryDimensions}). Thus,
our Theorem \ref{dim=G=L} generalizes Proposition 4.2 from
\cite{Lia1}.

In Section \ref{SectionCorollaries} we also generalize Theorem 5.2
from \cite{Lia1} and show that $K_0$ of $R$ is isomorphic to $K_0$
of the regular ring $Q$ of $R$. Specifically, in Corollary
\ref{moja K 0 theorem}, we show that the map $\mu:
\mathrm{Proj}(R)$ $\rightarrow$ $\mathrm{Proj}(Q)$ given by
$[P]\mapsto [P\otimes_{R}Q]$ is the isomorphism of monoids with
the inverse $[S]\mapsto [S\cap R^n]$ if $S$ is a direct summand of
$Q^n,$ and that $\mu$ induces the isomorphism $K_0(R)\cong
K_0(Q).$

\section{Class $\ce$ of Baer *-Rings}
\label{SectionClassCe}

\subsection{Basics.}\label{BasicDefn}
Let $R$ be a ring. $R$ is a {\em *-ring} (or {\em ring with
involution}) if there is an operation $*: R\rightarrow R$ such
that
\[(x+y)^*=x^*+y^*,\;\;\;(xy)^*=y^*x^*,\;\;\;(x^*)^*=x\;\;\;\mbox{ for all }x,y\in R.\]
If $R$ is also an algebra over $k$ with involution $*,$ then $R$
is an {\em *-algebra} if $(ax)^*=a^*x^*$ for $a\in k,$ $x\in R.$

An element $p$ of a *-ring $R$ is called a {\em projection} if $p$
is idempotent ($p^2=p$) and self-adjoint ($p^*=p$). There is a
partial ordering on the set of projections of $R$ defined by
$p\leq q$ iff $p=pq.$ The following conditions are equivalent
$p\leq q$; $p=qp$; $pR\subseteq qR;$ $R p \subseteq Rq.$

There is an equivalence relation on the set of projections of a
*-ring $R$ defined by $p\sim q$ iff $w^*w=p$ and $ww^*=q$ for some
$w\in R.$ Such an element $w$ is called a {\em partial isometry.}

Define another relation on the set of projections of a *-ring $R:$
$p \preceq q$, if $p$ is equivalent to a subprojection of $q$
(i.e. $p\sim r\leq q$ for some projection $r$). The relation
$\preceq$ is reflexive and transitive.

A {\em Rickart *-ring} is a *-ring $R$ such that, for every $x\in
R,$ the right annihilator $\ann_r(x)=\{y\in R\; |\; xy=0 \}$ is
generated by a projection $p$
\[\ann_r(x)=p R.\]

The projection $p$ from the above definition is unique. Also, if
$R$ is a Rickart *-ring, the left annihilator of each element $x$
of $R$ is $Rq$ for some (unique) projection $q$ since
$\ann_l(x)=(\ann_r(x^*))^*.$

Every element $x$ of a Rickart *-ring $R$ determines a unique
projection $p$ such that $xp=x$ and $\ann_r(x)=\ann_r(p)=(1-p)R$
and a unique projection $q$ such that $qx=x$ and
$\ann_l(x)=\ann_l(q)=R(1-q).$ $p$ is called the {\em right
projection of $x$} and is denoted by RP$(x).$ $q$ is the {\em left
projection of $x$} and is denoted by LP$(x).$

The involution in every Rickart *-ring is {\em proper}: $x^*x=0$
implies $x=0$ (Proposition 2, p. 13, \cite{Be2}). From this
condition it easily follows that a Rickart *-ring is a nonsingular
ring (since in a proper *-ring $R$ $\ann_r(x)\cap x^* R=0$ for all
$x\in R$). Let us also recall that a *-ring is called {\em
$n$-proper} if $x_1^*x_1+x_2^*x_2+\ldots +x_n^*x_n=0$ imply
$x_1=x_2=\ldots =x_n=0.$ It is easy to see that a *-ring $R$ is
$n$-proper if and only if the ring of $n\times n$ matrices over
$R$ is proper.

The projections in a Rickart *-ring form a lattice (Proposition 7,
p. 14 in \cite{Be2}).

A {\em Rickart $C^*$-algebra} is a $C^*$-algebra (complete normed
complex algebra with involution such that $||a^*a||=||a||^2$) that
is also a Rickart *-ring.

A {\em Baer *-ring} is a *-ring $R$ such that, for every nonempty
subset $S$ of $R,$ the right annihilator $\ann_r(S)=\{y\in R\; |\;
xy=0$ for all $x\in S \}$ is generated by a projection $p$
\[\ann_r(S)=p R.\]
Since $\ann_l(S)=(\ann_r(S^*))^*,$ it follows that $\ann_l(S)=R q$
for some projection $q.$

A *-ring is Baer *-ring if and only if it is Rickart *-ring and
the lattice of projections is complete (Proposition 1, p. 20,
\cite{Be2}).

If $R$ is a Baer *-ring and $\{p_i | i\in I\}$ is a nonempty
family of projections in $R$, then \begin{equation}(\inf_{i\in I}
p_i) R = \bigcap_{i\in I} p_i R\label{InfSup}\end{equation} This
is an easy exercise (Exercise 1A in \cite{Be2}).

A $C^*$-algebra that is a Baer *-ring is called an {\em
$AW^*$-algebra}.

If $H$ is a Hilbert space and $\B(H)$ the algebra of bounded
operators on $H,$ then $\B(H)$ is an $AW^*$-algebra. If $A$ is a
$*$-subalgebra of ${\B}(H)$ such that $A= A''$ where $A'$ is the
commutant of $A,$ then $A$ is called a {\em von Neumann algebra.}

A von Neumann algebra is an $AW^*$-algebra (\cite{Be2},
Proposition 9). The converse is not true (see \cite{D}) namely
there is an $AW^*$-algebra that cannot be represented as a von
Neumann algebra on any Hilbert space.

\subsection{Dimension.}\label{finiteBaer}
We now focus our attention on a special class of Baer *-rings.
\begin{itemize}
\item[(A1)] A Baer *-ring $R$ is {\em finite} if $x^*x=1$ implies
$xx^*=1$ for all $x\in R.$
\end{itemize}

The Baer *-ring $R$ satisfies the {\em generalized comparability
(GC)} axiom if: for every two projections $p$ and $q,$ there is a
central projection $c$ such that
\[ cp\preceq cq\;\;\mbox{ and }\;\;(1-c)q\preceq
(1-c)p.\]

We are interested in finite Baer *-rings with (GC) because of the
dimension function that we can define on the set of all
projections. Let $R$ be a finite Baer $*$-ring with (GC). Let $Z$
denote the center of $R.$ The projection lattice $P(Z)$ of $Z$ is
a complete Boolean algebra and, as such, may be identified with
the Boolean algebra of closed-open subspaces of a Stonian space
$X$. The space $X$ can be viewed as the set of maximal ideals in
$P(Z);$ $p\in P(Z)$ can be identified with the closed-open subset
of $X$ that consist of all maximal ideals that exclude $p.$

The algebra $C(X)$ of continuous complex-valued functions on $X$
is a commutative $AW^*$-algebra. An element $p\in P(Z)$ can be
viewed as an element of $C(X)$ by identifying $p$ with the
characteristic function of the closed-open subset of $X$ to which
$p$ corresponds.

If $R$ is an $AW^*$-algebra, then $Z$ is the closed linear span of
$P(Z)$ and we may identify $Z$ with $C(X).$

For more details on this construction, see \cite{Be2}.

\begin{thm} If $R$ is a finite Baer $*$-ring that satisfies (GC),
then there exist unique a function $d: P(R)\rightarrow C(X)$ such
that
\begin{itemize}
\item[(D1)] $p\sim q$ implies $d(p)=d(q),$

\item[(D2)] $D(p)\geq 0,$

\item[(D3)] $d(c)=c$ for every $c\in P(Z),$

\item[(D4)] $pq=0$ implies $d(p+q)=d(p)+d(q).$
\end{itemize}
The function $D$ will be called the {\em dimension function}. It
satisfies the following properties:
\begin{itemize}
\item[(D5)] $0\leq d(p)\leq 1,$

\item[(D6)] $d(cp)=cd(p)$ for every $c\in P(Z),$

\item[(D7)] $d(p)=0$ iff $p=0,$

\item[(D8)] $p\sim q$ iff $d(p)=(q),$

\item[(D9)] $p\preceq q$ iff $d(p)\leq d(q),$

\item[(D10)] If $p_i$ is an increasingly directed family of
projections with supremum $p$, then $d(p)= \sup d(p_i),$

\item[(D11)] If $p_i$ is an orthogonal family of projections with
supremum $p$, then $d(p)= \sum d(p_i).$
\end{itemize}
\label{DimensionExists}
\end{thm}

Chapter 6 of \cite{Be2} is devoted to the proof of this theorem.

\subsection{The Regular Ring of $R$}\label{Q}
Next, we would like to be able to enlarge our Baer *-ring to a
regular Baer *-ring. Recall that a ring $Q$ is {\em regular} if,
for every $x\in Q$ there is $y\in Q$ such that $xyx=x.$
Equivalently, a ring is regular if every right (left) module is
flat. A regular ring can also be characterized by the condition
that all finitely presented right modules are projective.

If $Q$ is a regular Rickart *-ring and $x\in Q,$ then $xQ=pQ$ for
some projection $p\in Q$ (Proposition 3, p. 229 in \cite{Be2}).

Every finite von Neumann algebra $A$ can be enlarged (in a
canonical way) to a regular ring $Q$ of certain unbounded
operators affiliated (densely defined and closed) with $A.$ In
chapter 8 of \cite{Be2}, this construction is generalized for a
certain class of finite Baer *-rings. The conditions that we must
impose onto a finite Baer *-ring $R$ in order to be able to follow
this construction are the following:
\begin{itemize}
\item[(A2)] $R$ satisfies {\em existence of projections
(EP)-axiom:} for every $0\neq x\in R,$ there exist an self-adjoint
$y\in\{x^*x\}''$ such that $(x^*x)y^2$ is a nonzero projection;

$R$ satisfies the {\em unique positive square root (UPSR)-axiom:}
for every $x\in R$ such that $x=x_1^*x_1+x_2^*x_2+\ldots
+x_n^*x_n$ for some $n$ and some $x_1,x_2,\ldots,x_n\in R$ (such
$x$ is called {\em positive}), there is a unique $y\in\{x^*x\}''$
such that $y^2=x$ and $y$ positive. Such $y$ is denoted by
$x^{1/2}.$

\item[(A3)] Partial isometries are addable.

\item[(A4)] $R$ is {\em symmetric}: for all $x\in R,$ $1+x^*x$ is
invertible.

\item[(A5)] There is a central element $i\in R$ such that $i^2=-1$
and $i^*=-i.$

\item[(A6)] $R$ satisfies the {\em unitary spectral (US)-axiom:}
for each unitary $u\in R$ such that RP$(1-u)=1,$ there exist an
increasingly directed sequence of projections $p_n\in\{u\}''$ with
supremum 1 such that $(1-u)p_n$ is invertible in $p_n Rp_n$ for
every $n.$

\item[(A7)] $R$ satisfies the {\em positive sum (PS)-axiom;} if
$p_n$ is orthogonal sequence of projections with supremum 1 and
$a_n\in R$ such that $0\leq a_n\leq f_n,$ then there is $a\in R$
such that $a p_n=a_n$ for all $n.$
\end{itemize}

By Theorem 1, p. 80, from \cite{Be2}, the generalized
comparability (GC) follows from (A2). Thus, all the Baer *-rings
satisfying (A1) and (A2) have the dimension function.

In the presence of (A2), the notion of positivity can be
simplified so that $x\in R$ is positive if and only if $x=y^*y$
for some $y\in R.$

(A2) -- (A5) imply that $R$ is $n$-proper (Lemma, p. 227 in
\cite{Be2}).

\begin{thm} If $R$ is a Baer *-ring satisfying (A1)-- (A7), then
there is a regular Baer *-ring $Q$ satisfying (A1) -- (A7) such
that $R$ is *-isomorphic to a *-subring of $Q$, all projections,
unitaries and partial isometries of $Q$ are in $R,$ and $Q$ is
unique up to *-isomorphism.

If $R$ is also an algebra over involutive field $F$, then so is
$Q.$ \label{RegularRing}
\end{thm}

This result is contained in Theorem 1, p. 217, Theorem 1 and
Corollary 1 p. 220, Corollary 1, p. 221, Theorem 1 and Corollary 1
p. 223, Proposition 3 p. 235, Theorem 1 p. 241, Exercise 4A p. 247
in \cite{Be2}.

A ring $Q$ as in Theorem \ref{RegularRing} is called the {\em
regular ring of Baer *-ring $R$.}

\begin{prop} If $R$ is a Baer *-ring satisfying (A1)-- (A7) with $Q$ its regular ring, then
\begin{enumerate}
\item $Q$ is the classical ring of quotients $Q_{\mathrm{cl}}(R)$
of $R.$

\item $Q$ is the maximal ring of quotients $Q_{\mathrm{max}}(R)$
of $R$ and, thus, self-injective and equal to the injective
envelope $E(R)$ of $R.$
\end{enumerate}
\label{Q=max=classical}
\end{prop}
\begin{pf}
1. First, let us show that $x\in R$ is a non-zerodivisor if and
only if it is invertible in $Q.$ It is easy to see that $x\in R$
that is invertible in $Q$ cannot have nontrivial left and right
annihilators. Conversely, if $x$ does not have a right inverse,
then the right annihilator of $x$ in $Q$ is nontrivial. Since $Q$
is Rickart, there is a nontrivial projection $p$ that generates
the right annihilator. But $p$ is in $R$ by Theorem
\ref{RegularRing}. Thus, $xp=0.$ The proof for the left handed
version is similar.

By Proposition 5, p. 241 of \cite{Be2}, for every $x\in Q$ there
is a partial isometry $w\in R$ such that $x=w (x^*x)^{1/2}.$
$(x^*x)^{1/2}$ is positive and thus self-adjoint. By Proposition
2, p. 228 \cite{Be2}, there is a unitary $u$ such that $1-u$ is
invertible in $Q$ and $(x^*x)^{1/2}=i(1+u)(1-u)^{-1}.$ But $u$ is
in $R$ by Theorem \ref{RegularRing}, so $1-u\in R.$ Thus, every
element $x\in Q$ can be represented as the right fraction $x=a
t^{-1},$ $a=w i (1+u)\in R,$ $t=1-u\in R.$

This proves the right Ore condition for $R$. Applying involution,
we have the left Ore condition. Since the non-zerodivisors of $R$
are invertible in $Q$, the isomorphism of $Q$ and
$Q_{\mathrm{cl}}(R)$ exists because of the universal property of
$Q_{\mathrm{cl}}(R).$ It is easy to check that the isomorphism is
a *-isomorphism.

2. In \cite{Pyle}, a *-extension of a finite Baer *-ring is
constructed that is (under suitable assumptions) *-isomorphic to
the maximal ring of quotients (\cite{Pyle}, Theorem 5.2) and
*-isomorphic to the regular ring of that finite Baer *-ring
(\cite{Pyle}, Theorem 5.3) with *-isomorphisms that fix the
original ring. Thus, to prove part 2. it is sufficient to show
that $R$ satisfies all the assumptions of Theorems 5.2 and 5.3
from \cite{Pyle}.

The assumptions for Theorem 5.2 are that the Baer *-ring is finite
(given by (A1)), every nonzero right ideal contains a nonzero
projection (guaranteed by EP, thus (A2)), LP $\sim$ RP (which
follows from (A2) by Corollary, p. 131 \cite{Be2}) and certain
condition called Utumi's condition. By Corollary 3.7 from
\cite{Pyle}, Utumi's condition is satisfied for every Baer *-ring
that is finite (A1), 2-proper (we have shown that (A2)-(A5) imply
$n$-proper for any positive $n$) and that (EP) and (SR) hold. (SR)
is an axiom that follows from (A2) and (A3) (see Exercise 7C p.
131, \cite{Be2}). Thus, all the assumptions are satisfied by $R.$

The assumptions for Theorem 5.3 are the same as (A1) -- (A6) with
the exception that (UPSR) in (A2) is replaced by (SR). But (SR)
follows from (A2) and (A3) and thus the assumptions of Theorem 5.3
hold for $R.$ \qed
\end{pf}

\subsection{Matrix Rings over $R$}\label{Matrices} Let $M_n(R)$ denotes the ring of
$n\times n$ matrices over $R$.

If $R$ is a Baer *-ring, the lattice of projections of $R$ is
complete. In order to ensure the completeness of lattice of
projections of $M_n(R)$ it is necessary for $M_n(R)$ to be Baer.
To ensure that we need two more axioms.

\begin{itemize}
\item[(A8)] $M_n(R)$ satisfies the {\em parallelogram law (P):}
for every two projections $p$ and $q,$ \[ p - \inf\{ p,q\}\sim
\sup\{ p,q\}-q.\]

\item[(A9)] Every sequence of orthogonal projections in $M_n(R)$
has a supremum.
\end{itemize}

If $Q$ is a regular ring of Baer *-ring $R$ that satisfies (A1)--
(A9), then $M_n(Q)$ is a regular Rickart *-ring that has the same
projections, unitaries and partial isometries as $M_n(R)$
(Propositions 2 and 3, p. 250 in \cite{Be2}). Thus, $Q$ satisfies
(A1) -- (A7) by Theorem \ref{RegularRing} and statements (A8) and
(A9) are true in $M_n(Q)$ (they are statements about projections
and the projections in $M_n(Q)$ and $M_n(R)$ are the same). Thus,
$Q$ satisfies (A1)-- (A9).

\begin{thm} If $R$ is a Baer *-ring satisfying (A1)-- (A9), then
$M_n(R)$ is a finite Baer *-ring with (GC).

If $Q$ is a regular Baer *-ring satisfying (A1)-- (A9), then
$M_n(Q)$ is a regular Baer *-ring.

\label{A1A9}
\end{thm}

This result is Theorem 1 and Corollary 2, p. 262 in \cite{Be2}.

\begin{cor}
\begin{enumerate}
\item If $R$ is a Baer *-ring satisfying (A1) -- (A7), then
$M_n(R)$ is semihereditary (i.e., every finitely generated
submodule of a projective module is projective or, equivalently,
every finitely generated ideal is projective) for every positive
$n.$

\item If $R$ is a Baer *-ring satisfying (A1) -- (A9), then the
lattice of projections of $M_n(R)$ is complete for every positive
$n.$
\end{enumerate}
\end{cor}
\begin{pf}
1. (A1) -- (A7) guarantees that $M_n(R)$ is a Rickart *-ring
(Theorem 1, p. 251 in \cite{Be2}).

A ring is right semihereditary if and only if the algebra of
$n\times n$ matrices is right Rickart for every positive $n$ (see
e.g. Proposition 7.63 in \cite{Lam}). Note that this result has a
corollary that if $R$ is right semihereditary, then $M_n(R)$ is
right semihereditary for every positive $n$ (simply identify
$M_m(M_n(R))$ with $M_{mn}(R)$ and use the result).

Thus, $M_n(R)$ is semihereditary for every positive $n.$

2. Since every Baer *-ring has a complete lattice of projections,
this is a simple corollary of the fact that $M_n(R)$ is a Baer
*-ring. \qed
\end{pf}

\begin{defn} Let $\ce$ be the class of Baer *-rings that satisfy
the axioms (A1) -- (A9.)
\end{defn}

Every finite $AW^*$-algebra satisfies the axioms (A1) -- (A9)
(remark 1, p. 249 in \cite{Be2}). Thus, the class $\ce$ contains
the class of all finite $AW^*$-algebras and, in particular, all
finite von Neumann algebras.

\section{Torsion Theories}
\label{SectionTT}

To study the properties of the class $\ce,$ we shall use a notion
that will facilitate the understanding of modules over a ring from
$\ce.$

We begin with a general setting: Let $R$ be any ring. A {\em
torsion theory} for $R$ is a pair $\tau = (\te, \ef)$ of classes
of $R$-modules such that
\begin{itemize}
\item[i)] $\homo_R(T,F)=0,$ for all $T \in \te$ and $F \in \ef.$
\item[ii)] $\te$ and $\ef$ are maximal classes having the property
$i).$
\end{itemize}
Thus, if $(\te, \ef)$ is a torsion theory, the class $\te$ is
closed under quotients, direct sums and extensions and the class
$\ef$ is closed under taking submodules, direct products and
extensions.

Conversely, if $\mem$ is a class of $R$-modules closed under
quotients, direct sums and extensions, then it is a torsion class
for a torsion theory $(\mem, \ef)$ where $\ef =
\{\;F\;|\;\homo_R(M,F)=0,\mbox{ for all }M\in\mem\;\}.$ Dually, if
$\mem$ is a class of $R$-modules closed under submodules, direct
products and extensions, then it is a torsion-free class for a
torsion theory $(\te, \mem)$ where $\te =
\{\;T\;|\;\homo_R(T,M)=0,\mbox{ for all }M\in\mem\;\}.$

The modules in $\te$ are called {\em $\tau$-torsion modules} (or
torsion modules for $\tau$) and the modules in $\ef$ are called
{\em $\tau$-torsion-free modules} (or torsion-free modules for
$\tau$).

If $\tau_1 = (\te_1, \ef_1)$ and $\tau_2 = (\te_2, \ef_2)$ are two
torsion theories, we say that $\tau_1$ is {\em smaller} than
$\tau_2$ ($\tau_1\leq\tau_2$) iff $\te_1\subseteq\te_2$
(equivalently $\ef_1\supseteq\ef_2$).

If $\mem$ is a class of $R$-modules, then torsion theory {\em
generated} by $\mem$ is the smallest torsion theory $(\te, \ef)$
such that $\mem\subseteq\te.$ The torsion theory {\em cogenerated}
by $\mem$ is the largest torsion theory $(\te, \ef)$ such that
$\mem\subseteq\ef.$

If $(\te, \ef)$ is a torsion theory for a ring $R$ and $M$ is a
$R$-module, there exists submodule $N$ such that $N\in\te$ and
$M/N\in\ef$ (Proposition 1.1.4 in \cite{Bland}). From this it
follows that every module $M$ has the largest submodule that
belongs to $\te$ (i.e. submodule generated by the union of all
torsion submodules of $M$). We call it the {\em torsion submodule}
of $M$ and denote it with $\te M$. The quotient $M/\te M$ is
called the {\em torsion-free quotient} and we denote it $\ef M.$

We say that a torsion theory $\tau = (\te, \ef)$ is {\em
hereditary} if the class $\te$ is closed under taking submodules.
A torsion theory is hereditary if and only if the torsion-free
class is closed under formation of injective envelopes
(Proposition 1.1.6, \cite{Bland}). Also, a torsion theory
cogenerated by a class of injective modules is hereditary (easy to
see) and, conversely, every hereditary torsion theory is
cogenerated by a class of injective modules (Proposition 1.1.17,
\cite{Bland}).

The notion of the closure of a submodule in a module is another
natural notion that can be related to a torsion theory. Let $M$ be
an $R$-module and $K$ a submodule of $M.$ The {\em closure}
$\cl_{\te}^M(K)$ of $K$ in $M$ with respect to the torsion theory
$(\te, \ef)$ is
\[\cl_{\te}^M(K) = \pi^{-1}(\te(M/K))\mbox{ where } \pi\mbox{ is
the natural projection }M\twoheadrightarrow M/K.\]

If it is clear in which module we are closing the submodule $K,$
we suppress the superscript $M$ from $\cl_{\te}^M(K)$ and write
just $\cl_{\te}(K)$. If $K$ is equal to its closure in $M,$ we say
that $K$ is {\em closed} submodule of $M$.

For more details on closure see Proposition 3.2 in \cite{Lia1}.

\subsection{Examples.}
\label{Examples}

\subsubsection{ } The torsion theory cogenerated by the injective
envelope $E(R)$ of $R$ is called the {\em Lambek torsion theory}.
We denote it $\tau_L.$ It is hereditary, as it is cogenerated by
an injective module, and faithful. Moreover, it is the largest
hereditary faithful torsion theory.

\subsubsection{ } The class of nonsingular modules over a ring $R$ is closed
under submodules, extensions, products and injective envelopes.
Thus, it is a torsion-free class of a hereditary torsion theory.
This theory is called the {\em Goldie torsion theory}. Let us
denote it with $\tau_G=(\T, \bigP).$

The Lambek theory is smaller than the Goldie theory (see example
3, p. 26 in \cite{Bland}). If $R$ is nonsingular, then the Lambek
and Goldie theories coincide (also see \cite{Bland} for details).

Here we mention a few results that we shall be using in the
sequel. By Corollary 7.30 in \cite{Lam}, if $M$ an $R$-module and
$K$ a submodule of $M,$ then the Goldie closure of $K$ in $M$ is
complemented in $M$. By Proposition 7.44 in \cite{Lam}, if $R$ is
a nonsingular ring and $M$ nonsingular $R$-module, then the Goldie
closure of $K$ in $M$ is the largest submodule of $M$ in which $K$
is essential. From this it follows that a submodule $K$ is Goldie
closed in $M$ if and only if $K$ is a complement in $M$.

The above has the following  result of R.E. Johnson (introduced in
\cite{John}) as a corollary.

\begin{cor}
Let $R$ be any ring and $M$ a nonsingular $R$-module. There is an
one-to-one correspondence \[\{\mbox{complements in }M\}
\longleftrightarrow\{\mbox{direct summands of }E(M)\}\] given by
$K \mapsto$ the Goldie closure of $K$ in $E(M)$ that is equal to a
copy of $E(K).$ The inverse map is given by $L\mapsto L\cap M.$
\label{Johnson}
\end{cor}

The proof can be found also in \cite{Lam} (Corollary 7.44').

\subsubsection{ } If $R$ is an Ore ring with the set of regular elements $T$ (i.e.,
$Tr \cap Rt \neq 0,$ for every $t \in T$ and $r\in R$), we can
define a hereditary torsion theory by the condition that a right
$R$-module $M$ is a torsion module iff for every $m\in M$, there
is a nonzero $t\in T$ such that $mt =0.$ This torsion theory is
called the {\em classical torsion theory of an Ore ring}. It is
faithful and so it is contained in the Lambek torsion theory.

\subsubsection{ } The class of flat modules is closed under extensions. If $R$
is semihereditary, the class of flat modules is closed under
direct products (Theorem 4.47 and Example 4.46 b), \cite{Lam}). If
$R$ is subflat (i.e. every submodule of a flat module is flat), it
is closed under submodules. Since every semihereditary ring $R$ is
subflat (Theorem 4.67, \cite{Lam}), semihereditary $R$ has a
torsion theory in which the class of all flat modules is the
torsion-free class. Denote this torsion theory with
$\tau_{\mathrm{flat}}.$

\subsubsection{ } Let $R$ be a subring of a ring $S$. Let us look at a
collection of all $R$-modules $M$ such that $S\otimes_R M = 0.$
This collection is closed under quotients, extensions and direct
sums. Moreover, if $S$ is flat as an $R$-module, then this
collection is closed under submodules and, hence, defines a
hereditary torsion theory. In this case we denote this torsion
theory by $\tau_S.$

From the definition of $\tau_S$ it follows that
\begin{itemize}
\item[1.] The torsion submodule of $M$ in $\tau_S$ is the kernel
of the natural map $M\rightarrow S \otimes_R M.$

\item[2.] All flat modules are $\tau_S$-torsion-free.
\end{itemize}
By 2., $\tau_S$ is faithful. Thus, $\tau_S$ is contained in the
Lambek torsion theory.

If $R$ is an Ore ring, then $\tau_{Q_{\mathrm{cl}}(R)}$ is the
classical torsion theory.

If $R$ is right semihereditary ring $R,$ the ring of maximal right
quotients that is left $R$-flat (Theorem 2.10 in
\cite{Sandomiersky}) and all torsion-free modules in
$\tau_{Q^r_{\mathrm{max}}(R)}$ are flat (Theorem 2.1 in
\cite{Turnidge}). Thus, if $R$ is Ore and semihereditary ring with
$Q_{\mathrm{cl}}(R)=Q_{\mathrm{max}}(R)$ (as is the case with any
ring from the class $\ce$) , then
\[\mbox{ Classical torsion theory }=\tau_{Q_{\mathrm{cl}}(R)}=\tau_{Q_{\mathrm{max}}(R)}=\tau_{\mathrm{flat}}.\]
In this case, let us denote this torsion theory by $(\smallt,
\p).$

\subsubsection{ } If $R$ is any ring, let $(\bnd, \unb)$ be the torsion theory cogenerated by the
ring $R$. We call a module in $\bnd$ a {\em bounded module} and a
module in $\unb$ an {\em unbounded module}. This theory is not
necessarily hereditary.

The Lambek and $(\bnd, \unb)$ torsion theory are related such that
$M$ is a Lambek torsion module if and only if every submodule of
$M$ is bounded. This is a direct corollary of the fact that
$\homo_R(M, E(R)) = 0$ if and only if $\homo_R(N, R)=0,$ for all
submodules $N$ of $M,$ that is an exercise in \cite{Cohn1}. Also,
it is easy to show that $(\bnd, \unb)$ is equal to the Lambek
torsion theory if and only if $(\bnd, \unb)$ is hereditary.

$(\bnd, \unb)$ is the largest torsion theory in which $R$ is
torsion-free. Thus, for a ring from class $\ce$
\[(\smallt,\p)\leq (\T, \bigP)= \mbox{ Lambek } \leq (\bnd, \unb).\]

\section{Torsion Theories for Rings from Class $\ce$}
\label{SectionTTforCe}

For the remainder of this section, let $R$ denote a ring from
class $\ce$ with $Q$ the regular ring of $R$. If $p$ is a matrix
from $M_n(R),$ we will identify $p$ with the $R$-map
$R^n\rightarrow R^n$ defined by $r\mapsto p r.$

\subsection{Splitting of $(\bnd, \unb)$ for Finitely Generated Modules.}

First, we shall show that $M=\bnd M\oplus \unb M$ for every
finitely generated $M$ and that $\unb M$ is finitely generated
projective. We need a few preliminary results.

\begin{lem} Let $P$ be a right $R$-module.
\begin{enumerate}
\item If $P$ is a submodule of $R^n,$ then the following
conditions are equivalent
\begin{itemize}
\item[i)] $P$ is a complement in $R^n.$

\item[ii)] There is a projection $p\in M_n(R)$ such that $P=\im
p.$

\item[iii)] $P$ is a direct summand of $R^n.$
\end{itemize}

\item $P$ is finitely generated projective if and only if there is
a nonnegative integer $n$ and a projection $p\in M_n(R)$ such that
$P=\im p.$
\end{enumerate}
\label{complements=summands}
\end{lem}
\begin{pf}
(1) i)$\Rightarrow$ ii) Let $P$ be a complement in $R^n.$ By
Corollary \ref{Johnson}, $E(P)$ is a direct summand of
$E(R^n)=E(R)^n=Q^n.$ The projection from $Q^n$ onto $E(P)$ is an
idempotent element $q\in M_n(Q)$ such that $\im q=E(P).$ Since
$M_n(Q)$ is a Rickart *-ring (by Theorem \ref{A1A9}), there is a
projection $p\in M_n(Q)$ such that $p M_n(Q)=\ann_r(1-q)=q
M_n(Q).$ But the projections in $M_n(Q)$ and $M_n(R)$ are the same
so $p\in M_n(R).$

$P$ is a complement, so $P=E(P)\cap R^n$ by Corollary
\ref{Johnson}. Thus, $P=p(Q^n)\cap R^n.$ Since $p\in M_n(R),$
$p(R^n)\subseteq R^n$ and so $p(R^n)\subseteq p(Q^n)\cap R^n=P.$
Conversely, if $p(r)$ is an element of $p(Q^n)\cap R^n=P,$ then
$p(r)\in R^n$ has unique decomposition as $p(r')+(1-p)(r'').$ But
that decomposition still holds in $Q^n.$ Thus, $p(r)=p(r')$ and
$(1-p)(r'')=0.$ Since $r'\in R^n,$ $p(r)=p(r')$ is in $p(R^n).$
This proves that $P=p(R^n).$

ii)$\Rightarrow$ iii) Trivial.

iii)$\Rightarrow$ i) Trivial.

(2) If $P$ is finitely generated projective, then there is a
nonnegative integer $n$ such that $P$ is a direct summand in
$R^n.$ Then $P=\im p$ for some projection $p\in M_n(R)$ by (1).
The converse is obvious. \qed
\end{pf}

The following lemma asserts that we can separate a direct summand
and an element in the image of a projection out of the direct
summand, with an $R$-valued map. This will turn out to be the key
ingredient in the proof that a finitely generated module $M$
splits as $\bnd M\oplus \unb M.$

\begin{lem} If $P$ is a direct summand of $R^n,$ $p\in M_n(R)$ a
projection, and $a\in R^n$ any element such that $p(a)\notin P$,
then there is a map $f\in \homo_R(R^n, R)$ such that $f(P)\equiv
0$ and $f(p(a))\neq 0.$ \label{separation}
\end{lem}
\begin{pf} Let $S$ be the complement of $P$ and $pr_S$ be the
projection of $R^n$ onto $S$. $p(a)=r_P+r_S$ where $r_P\in P$ and
$r_S\in S.$ Since $p(a)\notin P,$ $r_S$ is nontrivial. Let $(q_1,
q_2, \ldots, q_n)$ be the coordinates of $r_S$ in the standard
basis. Define the map $q: R^n\rightarrow R$ by
\[g: (a_1, a_2,\ldots, a_n)\mapsto \sum_{i=1}^ n q^*_i a_i \]

Now define map $f\in\homo_R(R^n, R)$ as $f=g\circ pr_{S}.$ Clearly
$P\in \ker f.$ $f(p(a))=q(r_S)=\sum_{i=1}^ n q^*_i q_i\neq 0$
since $R$ is $n$-proper and $r_S\neq 0.$ \qed
\end{pf}

The idea here is to study finitely generated projective modules by
treating the projections and benefit from the nice properties of
the matrix rings over $R.$ This is possible since if $p,q\in
M_n(R),$ then \begin{equation} p M_n(R)=q M_n(R)\;\;\mbox{ if and
only if }\;\;p(R^n)=q(R^n)\label{ImagesOfProjections}
\end{equation}
as basic matrix algebra shows. Now, we can understand the closures
in $R^n$ better.

\begin{prop} If $P$ is a submodule of $R^n,$
then the following sets are equal and are direct summands of
$R^n$:
\begin{enumerate}
\item $\cl_{\bnd}(P)=\{x\in R^n | f(x)=0$ for all
$f\in\homo_R(R^n,R)$ s.t. $P\subseteq \ker f\}=\bigcap\{\ker f |
f\in\homo_R(R^n,R)$ s.t. $P\subseteq \ker f\},$

\item $\bigcap\{ S | S$ is a direct summand of $R^n$ and
$P\subseteq S\},$

\item $\inf\{p | p\in M_n(R)$ a projection with $P\subseteq p(R^n)
\} (R^n) = \bigcap\{ p(R^n) | p\in M_n(R)$ a projection with
$P\subseteq p(R^n) \},$

\item $\cl_{\T}(P)  =$ (largest submodule of $R^n$ in which $P$ is
essential) = (smallest submodule of $R^n$ that contains $P$ and
that is a complement in $R^n) = E(P)\cap R^n$.

Moreover, if $P$ is a right ideal in $R$, then the above sets are
equal to $\ann_r(\ann_l(P)).$
\end{enumerate}
\label{VariousClosures}
\end{prop}
\begin{pf}
All sets in (1) are equal by the definition of closure in the
torsion theory $(\bnd, \unb).$ Also, if $P$ is a right ideal of
$R,$ then it is easy to see that
$\cl_{\bnd}(P)=\ann_r(\ann_l(P)).$

The sets in (3) are equal by formula (\ref{InfSup}) in subsection
\ref{BasicDefn} and formula (\ref{ImagesOfProjections}) above.

The first three sets in (4) are equal by Proposition 7.44 from
\cite{Lam}. From Corollary \ref{Johnson} also follows that
$\cl_{\T}(P)=E(\cl_{\T}(P))\cap R^n.$ But $P\subseteq_e
\cl_{\T}(P)$ and so $E(P)=E(\cl_{\T}(P)).$ Thus
$\cl_{\T}(P)=E(\cl_{\T}(P))\cap R^n=E(P)\cap R^n.$

$(2) \subseteq (1)$ since every $f$ as in (1) determines one $S$
as in (2).

$(2) = (3)$ by Lemma \ref{complements=summands}.

$(2) = (4)$ The intersection of complements is a complement by
Proposition 7.44 in \cite{Lam}. Thus, the set in (4) is the
intersection of all complements in $R^n$ containing $P$. But by
Lemma \ref{complements=summands}, this is the same as the
intersection of all direct summands of $R^n$ containing $P$.
Moreover, the set in (4) is a complement itself and, thus a direct
summand in $R^n.$

The set in (1) is also a complement since it is the intersection
of complements. Thus, the set in (1) is a direct summand as well.

Let us show now that $(1) \subseteq (2).$ Since the set in (1) is
a direct summand, there is a projection $p\in M_n(R)$ such that
$\cl_{\bnd}(P)=p(R^n).$ To show $(1) \subseteq (2)$ it is
sufficient to show that $p(a)$ is contained in all direct summands
of $R^n$ that contain $P$ for all $a\in R^n$ (because then
$p(R^n)\subseteq (2)).$

Suppose the contrary: there is $a\in R^n$ and a direct summand $S$
such that $P\subseteq S$ and $p(a)\notin S.$ By Lemma
\ref{separation}, there is a map $f\in \homo_R(R^n, R)$ such that
$f(S)\equiv 0$ and $f(p(a))\neq 0.$ But $p(a)$ is in
$\cl_{\bnd}(P)$ and $f$ is a map such that $P\subseteq S\subseteq
\ker f$ so $p(a)$ has to be in the kernel of $f$ as well.
Contradiction. Thus, $p(a)\in S.$ \qed
\end{pf}

\begin{thm} If $M$ is finitely generated $R$-module and $K$ submodule of $M$,
then $M/\cl_{\bnd}(K)$ is finitely generated projective and
$\cl_{\bnd}(K)$ is a direct summand of $M.$ In particular, for
$K=0$ we have that $\unb M$ is finitely generated projective and
$M=\bnd M\oplus \unb M.$ \label{SplittingOfBnd}
\end{thm}
\begin{pf}
If $M$ is $R^n$, $\cl_{\bnd}(K)$ is a direct summand of $M$ by
Proposition \ref{VariousClosures}. Moreover, the inclusion of
$\cl_{\bnd}(K)$ in $M$ splits since $\cl_{\bnd}(K)=p(R^n)$ for
some projection $p\in M_n(R^n).$ Thus the claim follows for
$M=R^n$.

Now let $M$ be any finitely generated module. There is a
nonnegative integer $n$ and an epimorphism $f: R^n\rightarrow M.$

First, we shall show that $\cl_{\bnd}(f^{-1}(K)) =
f^{-1}(\cl_{\bnd}(K)).$

Let $x$ be in $\cl_{\bnd}(f^{-1}(K)).$ Then $g(x)=0,$ for every
$g\in\homo_{R}(R^n, R)$ such that $f^{-1}(K)\subseteq\ker g.$ We
need to show that $f(x)$ is in $\cl_{\bnd}(K),$ i.e. that
$h(f(x))=0$ for every $h\in\homo_{R}(M, R)$ with $K\subseteq\ker
h.$ Let $h$ be one such map. Letting $g = hf,$ we obtain a map in
$\homo_{R}(R^n, R)$ such that $g(f^{-1}(K))=hff^{-1}(K) = h(K)$
(since $f$ is onto). But $h(K)=0, $ and so $f^{-1}(K)\subseteq
\ker g.$ Hence, $g(x) = 0$ i.e. $h(f(x))=0.$

To show the converse, let $x$ be in $f^{-1}(\cl_{\bnd}(K)).$ Then
$h(f(x)) = 0$ for every $h$ $\in$ $\homo_{R}(M, R)$ such that
$K\subseteq\ker h.$ We need to show that $g(x) = 0$ for every
$g\in\homo_{R}(R^n, R)$ such that $f^{-1}(K)\subseteq\ker g.$ Let
$g$ be one such map. Since $f^{-1}(0)\subseteq
f^{-1}(K)\subseteq\ker g,$ we have $\ker f\subseteq\ker g.$ This
condition enables us to define a homomorphism $h: M\rightarrow R$
such that $h(f(p)) = g(p)$ for every $p\in R^n.$ Then $h(K) = h(
f(f^{-1}(K))) = g(f^{-1}(K)) = 0,$ and so $h(f(x)) = 0.$ But from
this $g(x) =0.$

It is easy to see that $f: R^n\rightarrow M$ induces an
isomorphism $R^n/f^{-1}(\cl_{\bnd}(K))$ $\rightarrow$
$M/\cl_{\bnd}(K).$ But $ \cl_{\bnd}(f^{-1}(K))
=f^{-1}(\cl_{\bnd}(K)),$ so we obtain that $M/\cl_{\bnd}(K)$ is
finitely generated projective (since $R^n/\cl_{\bnd}(f^{-1}(K))$
is). So $0\rightarrow\cl_{\bnd}(K)\rightarrow M\rightarrow
M/\cl_{\bnd}(K)\rightarrow 0 $ splits. \qed \end{pf}

\subsection{Dimension.}

Given that the dimension function on $R$ and on all rings $M_n(R)$
exist (Theorem \ref{DimensionExists} and Theorem \ref{A1A9}) it
would be desirable to have the dimensions on $M_m(R)$ and $M_n(R)$
agree for $m\geq n$ i.e. $d_m\upharpoonleft_{M_n(R)}=d_n$ for all
$m\geq n.$

The dimension on $R$ is determined by its values on the central
projections (see chapter 6 of \cite{Be2}). The centers of $R$ and
$M_n(R)$ are isomorphic under the identification of
diag$(a,a,...,a)\in Z(M_n(R))$ with $a\in Z(R).$ If we identify
diag$(a,a,...,a)\in Z(M_n(R))$ with $na\in Z(R),$ we get the
desired result on the dimensions.

Now let us define the function $\dim_R$ on the class of all right
$R$ modules $\mbox{Mod}_R$ and values in $C(X)$ by
\begin{enumerate}
\item If $P$ is a finitely generated projective $R$-module,  then
there is a nonnegative integer $n$ and a matrix $p\in M_n(R)$ such
that $p^2=p^*=p$ and $\im p\cong P$. It is clear that an
idempotent matrix $q$ with image isomorphic to $p$ exist. Choose
$p$ to be the projection such that $p M_n(R)=\ann_r(1-q).$ Recall
that we can do that because $M_n(R)$ is a Rickart *-ring. Then
define
\[\dim_R(P)=d(p).\]

The values of $\dim_R$ are in $C_{[0,\infty)}(X),$ the algebra of
functions from $C(X)$ with values in $[0,\infty).$ The algebra
$C_{[0,\infty)}(X)$ is a boundedly complete lattice with respect
to the pointwise ordering (see pages 161 and 162 in \cite{Be2}).
Note, however, that the infinite lattice operations might differ
from the pointwise operations.

\item If $M$ is any $R$-module, define
\[\dim_R'(M)=\sup \{ \dim_R(P)\; |\; P \mbox{
fin. gen. projective submodule of }M\}\] where the supremum on the
right side is an element of $C(X)$ if it exists and is a new
symbol $\infty$ otherwise. We define
$a+\infty=\infty+a=\infty=\infty+\infty$ and $a\leq \infty$ for
every $a\in C(X).$
\end{enumerate}

Our first goal is to show that the following theorem (proven by
Wolfgang L\"{u}ck in \cite{Lu2}) holds for $R$ with $[0,\infty)$
replaced by $C_{[0,\infty)}(X)$ and $[0,\infty]$ replaced by
$C_{[0,\infty)}(X)\cup\{\infty\}.$

\begin{thm} Let $R$ be a ring such that there exist a dimension
function $\dim$ that assigns to any finitely generated projective
right $R$-module an element of $[0, \infty)$ and such that the
following two conditions hold
\begin{itemize}
\item[(L1)] If $P$ and $Q$ are finitely generated projective
modules, then
\[P\cong Q \Rightarrow \dim(P)=\dim(Q)\]
\[\dim(P\oplus Q)=\dim(P)+\dim(Q),\]

\item[(L2)] If $K$ is a submodule of finitely generated projective
module $Q,$ then $\cl_{\bnd}(K)$ is a direct summand of $Q$ and
\[\dim(\cl_{\bnd}(K))=\sup\{\dim(P)\;|\; P\mbox { is a fin. gen. projective submodule of } K\}.\]
\end{itemize}
Then, for every $R$-module $M$, we can define a dimension
\[\dim_R'(M)=\sup \{ \dim_R(P)\; |\; P \mbox{ fin. gen.
projective submodule of }M\}\in[0,\infty]\] that satisfies the
following properties:
\begin{enumerate}
\item Extension: $\dim(P)=\dim'(P)$ for every finitely generated
projective module $P.$

\item Additivity: If $\;0\rightarrow M_0\rightarrow M_1\rightarrow
M_2\rightarrow 0$ is a short exact sequence of $R$-modules, then
\[ \dim'(M_1)= \dim'(M_0)+\dim'(M_2).\]

\item Cofinality: If $M = \bigcup_{i \in I}M_i$ is a directed
union, then \[\dim'(M) = \sup\{\;\dim'(M_i)\; |\; i\in I\;\}.\]

\item Continuity: If $K$ is a submodule of a finitely generated
module $M$, then \[\dim'(K)=\dim'(\cl_{\bnd}(K)).\]

\item If $M$ is a finitely generated module, then
\[\dim'(M)=\dim(\unb M)\;\;\mbox{ and }\;\;\dim'(\bnd M)=0.\]

\item The dimension $\dim'$ is uniquely determined by (1) -- (4).
\end{enumerate}
\label{LueckDimension}
\end{thm}

For proof see Theorem 6.7, p 239 of \cite{Lu_book} or Theorem 0.6
and Remark 2.14 in \cite{Lu2}.

First, we show that the condition (L1) from Theorem
\ref{LueckDimension} holds for $R.$

\begin{prop} If $P$ and $S$ are finitely generated projective
$R$-modules, then
\begin{enumerate}
\item $P\cong S\;\; \mbox{ if and only if
}\;\;\dim_R(P)=\dim_R(S),$

\item $\dim_R(P\oplus S)=\dim_R(P)+\dim_R(S).$
\end{enumerate}
\label{LuecksFirstCondition}
\end{prop}

\begin{pf}
(1) Let $P\cong S.$ Let $p$ and $s$ be projections such that
$\dim_R(P)=d(p)$ and $\dim_R(S)=d(s).$ $p$ and $s$ might be
matrices of different size. Then there is an integer $n$ such that
\[
p_n = \left[\begin{array}{cc} p & 0\\
0 & 0 \end{array}\right] \mbox{ and } s_n = \left[\begin{array}{cc} s & 0\\
0 & 0 \end{array}\right] \] are both in $M_n(R)$ and there is an
invertible matrix $u\in M_n(R)$ such that $u p_n=s_n u$ (see Lemma
1.2.1. in \cite{Rosenberg} for details).

Similar elements are algebraically equivalent (i.e $a$ is
algebraically equivalent to $b$ iff $xy=a,$ $yx=b$ for some $x$
and $y$). Algebraic equivalence implies $\sim$ equivalence in all
*-rings with (SR) (Exercise 8A, p. 9 in\cite{Be2}). Since (SR)
holds if (A2) and (A3) hold, we have that $p_n\sim s_n.$ Thus,
$d(p)=d(p_n)=d(s_n)=d(s).$

Conversely, if $\dim_R(P)=\dim_R(S),$ then
$d(p)=d(p_n)=d(s_n)=d(s)$ (we might have to enlarge $p$ and $s$
again). So $p_n\sim s_n.$ It is easy to see (Exercise 5A, p. 8 in
\cite{Be2}) that then $\im p_n$ is isomorphic to $\im s_n.$ But
then $P$ is isomorphic to $S.$

(2) Let $P$ and $S$ be finitely generated projective modules with
$p$ and $s$ projections such that $\dim_R(P)=d(p)$ and
$\dim_R(S)=d(s).$ Then we can use
\[p\oplus s = \left[\begin{array}{cc} p & 0\\
0 & s \end{array}\right] \] to compute the dimension of $P\oplus
S.$ There is an integer $n$ such that \[
p_n = \left[\begin{array}{cc} p & 0\\
0 & 0 \end{array}\right] \mbox{ and } s_n = \left[\begin{array}{cc} 0 & 0\\
0 & s \end{array}\right] \] are both in $M_n(R).$ Then, $p_n
s_n=s_n p_n=0$ and so $\dim_R(P\oplus S)=d(p\oplus
s)=d(p_n+s_n)=d(p_n)+d(s_n)=d(p)+d(s).$ \qed \end{pf}

Note that this Proposition implies that \[\dim_R(P)=0\;\;\mbox{
iff }\;\;P=0\] for every finitely generated projective module $P.$

In order to prove that $R$ satisfies condition (L2) from Theorem
\ref{LueckDimension}, we need two lemmas.

Recall that the regular ring $Q$ of $R$ is also in the class $\ce$
(see subsection \ref{Matrices}). Thus, we can define its dimension
function $\dim_Q.$ The following result relates the dimensions of
$R$ and $Q$ and is leading us one step closer to (L2) of Theorem
\ref{LueckDimension}.

\begin{lem}
\begin{enumerate}
\item If $P$ is a direct summand of $R^n,$ then
\[\dim_R(P)=\dim_Q(E(P))\]

\item If $S$ is a direct summand of $Q^n,$ then
\[\dim_Q(S)=\dim_R(S\cap R^n).\]

\item If $S$ is a submodule of $Q^n,$ then
\[\dim_Q(\cl_{\bnd}(S))=\sup\{d(q)\;|\; q\in M_n(R)\mbox{ a projection, }q(Q^n)\subseteq S\}\]

\item If $P$ is a submodule of $R^n,$ then
\[\dim_R(\cl_{\bnd}(P))=\sup\{d(q)\;|\; q\in M_n(R) \mbox{ a projection, }q(R^n)\subseteq P\}\]
\end{enumerate}
\label{ClosureAndSup}
\end{lem}
\begin{pf}
(1) If $P$ is a direct summand of $R^n,$ $P= p(R^n)$ for some
projection $p\in M_n(R)$ by Lemma \ref{complements=summands}. By
definition of $\dim_R$ it follows that $\dim_R(P)=d(p).$

From the proof of i)$\Rightarrow$ ii) in Lemma
\ref{complements=summands}, it follows that $p(R^n)=p(Q^n)\cap
R^n.$ Thus, $P=p(Q^n)\cap R^n,$ and so $E(P)=E(p(Q^n)\cap
R^n)=p(Q^n)$ by Corollary \ref{Johnson}. Thus $\dim_Q(E(P))=d(p).$

(2) If $S$ is a direct summand of $Q^n,$ $S= p(Q^n)$ for some
projection $p\in M_n(R).$ Then $\dim_Q(S)=d(p).$ Then, $S\cap
Q^n=p(Q^n)\cap R^n=p(R^n)$ again by the proof of i)$\Rightarrow$
ii) in Lemma \ref{complements=summands}. Thus, $\dim_R(S\cap
Q^n)=d(p).$

(3) First we shall show that
\[\cl_{\bnd}(S) = \sup\{ q | q\in M_n(R)\mbox{ a projection with }q(Q^n)\subseteq
S\}(Q^n).\]

Let $p$ denote the projection $\sup\{ q | q\in M_n(R)\mbox{ a
projection with }q(Q^n)\subseteq S\},$ and $r$ denote the
projection such that $\cl_{\bnd}(S)=r(Q^n).$ We shall show that
$p=r.$

Since $q(Q^n)\subseteq S\subseteq r(Q^n)$ for all projections $q$
with $q(Q^n)\subseteq S$, then $q\leq r$ for all such $q$ and so
$p\leq r.$

Conversely, it is sufficient to show $S\subseteq p(Q^n)$ since
then $p(Q^n)\supseteq\inf\{q | q\in M_n(R)$ a projection with
$S\subseteq q(Q^n) \}(Q^n)=\cl_{\bnd}(S)=r(Q^n)$ by Proposition
\ref{VariousClosures} and thus $p\geq r.$ So, let $x\in S.$
Consider a matrix $X\in M_n(Q)$ such that the entries in the first
column are coordinates of $x$ in the standard basis and the
entries in all the other columns equal zero. Since $M_n(Q)$ is a
regular Rickart *-ring (Theorem \ref{A1A9}), there is a projection
$p_x\in M_n(Q)$ such that $X M_n(Q)= p_x M_n(Q).$ But $x\in
X(Q^n)=x Q\subseteq S$ and so we have that $p_x(Q^n)\subseteq S$
for all $x\in S.$ So, $p_x\leq p$ for all $x\in S.$ Thus, $x\in
p_x(Q^n)\subseteq p(Q^n)$ for all $x\in S$ and so $S\subseteq
p(Q^n).$

Now it is easy to see that
\[\begin{array}{rcl}\dim_Q(\cl_{\bnd}(S)) & = & \dim_Q(\sup\{ q | q\in
M_n(R)\mbox{ a projection with }q(Q^n)\subseteq S\}(Q^n))\\
& = & \sup\{ d(q) | q\in M_n(R)\mbox{ a projection with
}q(Q^n)\subseteq S\}\end{array}\] by property (D10) of Theorem
\ref{DimensionExists}.

(4) Let $p$ denote the projection $\sup\{ q | q\in M_n(R)\mbox{ a
projection with }q(R^n)\subseteq P\},$ and $r$ denote the
projection such that $\cl_{\bnd}(P)=r(R^n).$ Since
$q(R^n)\subseteq P\subseteq r(R^n),$ for all projections $q$ such
that $q(R^n)\subseteq P,$ then $q\leq r.$ Thus, $p\leq r$ and so
\[\dim_R(\cl_{\bnd}(P)) = d(r) \geq d(p)
 = \sup\{ d(q) | q\in M_n(R)\mbox{ projection, }q(R^n)\subseteq P\} \]

For the converse, first note that
$E(\cl_{\bnd}(P))=E(\cl_{\T}(P))$ by Proposition
\ref{VariousClosures}, $E(\cl_{\T}(P))=E(P)$ since $P$ is
essential in $\cl_{\T}(P),$ and $E(P)=E(P)\cap
Q^n=\cl^{Q}_{\bnd}(P)$ again by Proposition \ref{VariousClosures}.
Thus, $E(\cl_{\bnd}(P))=\cl^{Q}_{\bnd}(P).$ Then,

\[\begin{array}{rcl}
\dim_R(\cl_{\bnd}(P)) & = & \dim_Q(E(\cl_{\bnd}(P)))\hskip1.6cm
(\mbox{by
part (1)})\\
 & = & \dim_Q(\cl^{Q}_{\bnd}(P))) \hskip2cm (\mbox{by remark above})\\
 & = & \sup\{ d(q) | q\in M_n(R)\mbox{ proj.,
}q(Q^n)\subseteq P\} \hskip.3cm (\mbox{by part (3)})\\
 & \leq & \sup\{ d(q) | q\in M_n(R)\mbox{ proj.,
}q(R^n)\subseteq P\} \hskip.3cm (q(R^n)\subseteq q(Q^n)).
\end{array}\]
\qed
\end{pf}

\begin{lem}
\begin{enumerate}
\item If $P$ is a finitely generated projective module in $R^n$,
then \[\dim_R(P)=\dim_R(\cl_{\bnd}(P)).\]

\item If $P$ and $S$ are finitely generated projective modules in
$R^n,$ then

\[P\subseteq S\mbox{ imples }\;\;\dim_R(P)\leq \dim_R(S).\]
\end{enumerate}
\label{MonotonyForFGP}
\end{lem}
\begin{pf}
(1) Let $p\in M_n(R)$ be a projection such that $p(R^n)\cong P$
and $q$ be a projection such that $q(R^n)=\cl_{\bnd}(P).$ To prove
the claim, it is sufficient to show that $p\sim q$ (since then
$\dim_R(P)=d(p)=d(q)=\dim_R(\cl_{\bnd}(P))$). For $p\sim q,$ it is
sufficient to show that $q(Q^n)\cong p(Q^n)$ by the same argument
we used in the proof of part (1) of Proposition
\ref{LuecksFirstCondition}. $q(Q^n)=E(P)$ since $q(Q^n)\cap
R^n=q(R^n)=\cl_{\bnd}(P)=E(P)\cap R^n$ (by Corollary
\ref{Johnson}). Thus,
\[q(Q^n)=E(P)\cong E(p(R^n))=E(p(Q^n)\cap R^n)=p(Q^n).\]

(2) Let $p$ be a projection such that $p(R^n)=\cl_{\bnd}(P)$ and
$s$ a projection with $s(R^n)=\cl_{\bnd}(S).$ $P\subseteq S$
implies $p(R^n)=\cl_{\bnd}(P)\subseteq \cl_{\bnd}(S)=s(R^n).$
Thus, $p\leq s$ and so $d(p)\leq d(s).$ Hence
\[\dim_R(P)=\dim_R(\cl_{\bnd}(P))=d(p)\leq
d(s)=\dim_R(\cl_{\bnd}(S))=\dim_R(S)\] by part (1). \qed \end{pf}

Now we can prove he formula from Condition (L2), Theorem
\ref{LueckDimension}.

\begin{prop} If $K$ is a submodule of a finitely generated projective
module $S,$ then
\[\dim_R(\cl^S_{\bnd}(K))=\sup\{\dim_R(P)\;|\; P\mbox { is a fin. gen. projective submodule of } K\}\]
\label{LueckSecondCondition}
\end{prop}
\begin{pf} Since $S$ is a finitely generated projective, there is a nonnegative integer $n$
such that $S$ is a direct summand of $R^n.$ $\cl^S_{\bnd}(K)$ is a
direct summand of $S$ by Theorem \ref{SplittingOfBnd}. Thus,
$\cl^{R^n}_{\bnd}(K)\subseteq \cl^S_{\bnd}(K)$ by Proposition
\ref{VariousClosures}. $S\subseteq R^n$ implies
$\cl^S_{\bnd}(K)\subseteq \cl^{R^n}_{\bnd}(K).$ Hence,
$\cl^S_{\bnd}(K) = \cl^{R^n}_{\bnd}(K)$ so we can work in $R^n.$

\[\begin{array}{rcl}
\dim_R(\cl^{R^n}_{\bnd}(K)) & = & \sup\{d(q)\;|\; q\mbox { a
projection in }M_n(R)\mbox{ such that }q(R^n)\subseteq K\}\\

& \leq & \sup\{\dim_R(P)\;|\; P\mbox { is a fin. gen. proj.
submodule of } K\}
\end{array}\]
by Lemma \ref{ClosureAndSup}.

Conversely,
\[\begin{array}{lcr}
\sup\{\dim_R(P)\;|\; P\mbox { is a fin. gen. projective submodule
of } K\} & \leq & \\

\sup\{\dim_R(P)\;|\; P\mbox { is a fin. gen. projective submodule
of } \cl^{R^n}_{\bnd}(K)\} & = & \\
\dim_R(\cl^{R^n}_{\bnd}(K)). & &
\end{array}\]

The last equality holds since we have monotony for dimensions of
finitely generated projective modules by the Lemma
\ref{MonotonyForFGP} (this gives us $\leq$). The converse follows
since $\cl^{R^n}_{\bnd}(K)$ is finitely generated projective by
Theorem \ref{SplittingOfBnd}. \qed \end{pf}

Finally, we can prove that our dimension is just as in L\"{u}ck's
Theorem \ref{LueckDimension}. Recall that we need to replace
$[0,\infty)$ by $C_{[0,\infty)}(X)$ and $[0,\infty]$ by
$C_{[0,\infty)}(X)\cup\{\infty\}$ in the formulation of the
theorem. Luckily, that will not influence the proof.

\begin{thm} Theorem
\ref{LueckDimension} holds for $R$ and $\dim_R.$
\label{DimensionInR}
\end{thm}
\begin{pf}
$R$ satisfies condition (L1) of Theorem \ref{LueckDimension} by
Proposition \ref{LuecksFirstCondition} and the condition (L2) by
Theorem \ref{SplittingOfBnd} and by Proposition
\ref{LueckSecondCondition}.

The Extension property holds by Proposition
\ref{LueckSecondCondition} and Lemma \ref{MonotonyForFGP}. The
proof of the rest of the theorem is identical to the original
proof by L\"{u}ck  (proof of Theorem 6.7, p 239 of \cite{Lu_book}
or Theorem 0.6 and Remark 2.14 in \cite{Lu2}). \qed \end{pf}

From now on we shall not distinguish between $\dim_R$ and
$\dim_R'$ and will denote them both by $\dim_R.$

Further, the dimension has the following properties.
\begin{cor}
If $M = \bigoplus_{i\in I} M_i,$ then $\dim_R(M) = \sum_{i \in
I}\dim_R(M_i).$ \label{FurtherPropOfDim}
\end{cor}
\begin{pf}
This is an easy corollary of Cofinality of $\dim_R.$ \qed
\end{pf}

Theorem \ref{DimensionInR} and Corollary \ref{FurtherPropOfDim}
enable us to define another torsion theory:  for an $R$-module $M$
define $\T_{\dim} M$ as the submodule generated by all submodules
of $M$ of zero dimension. It is zero-dimensional by Cofinality of
$\dim_R.$ So, $\T_{\dim} M$ is the largest submodule of $M$ of
zero dimension.

Let us denote the quotient $M/ \T_{\dim} M$ by $\bigP_{\dim} M$.

The class $ \T_{\dim} = \{M \in \mathrm{Mod}_{R} | M = \T_{\dim}
M\}$ is closed under submodules, quotients and extensions by
additivity of dimension. The closure under the formation of direct
sums follows from Corollary \ref{FurtherPropOfDim}. Thus,
$\T_{\dim}$ defines a hereditary torsion theory with torsion-free
class equal to $\bigP_{\dim} = \{M \in \mathrm{Mod}_{R} | M =
\bigP_{\dim} M\}.$

Since $R$ is semihereditary and a nontrivial finitely generated
projective module has nontrivial dimension, $R$ is in
$\bigP_{\dim}.$ Thus, the torsion theory $(\T_{\dim},
\bigP_{\dim})$ is faithful. Since the Lambek torsion theory is the
largest hereditary and faithful theory,
\[(\T_{\dim}, \bigP_{\dim})\leq (\T, \bigP)\leq (\bnd, \unb).\]

\begin{thm} \[(\T_{\dim}, \bigP_{\dim})= (\T, \bigP).\]

If $M$ is finitely generated, then $\T_{\dim}(M)=\bnd M.$ Thus,
$M$ splits as $\T_{\dim} M\oplus \bigP_{\dim} M.$ \label{dim=G=L}
\end{thm}
The theorem allows us to drop the subscript $\dim$ from
$(\T_{\dim}, \bigP_{\dim}).$ The proof is the same as the proof of
Proposition 4.2 from \cite{Lia1}. We quote it for completeness.
\begin{pf} If $M$ is finitely generated, then $\dim_R(\bnd
M)=\dim_R(\cl_{\bnd}(0))=0$ by Continuity Property of dimension.
Thus, $\bnd M\subseteq \T_{\dim} M.$ Since the converse always
holds, the claim follows. The splitting follows from Theorem
\ref{SplittingOfBnd}.

Since $(\T_{\dim}, \bigP_{\dim})\leq (\T, \bigP)=$ Lambek torsion
theory, to prove the equality it is sufficient to prove that every
Lambek torsion module $M$ has dimension zero. Recall that $M$ is
Lambek torsion module iff all submodules of $M$ are bounded. This
means that all finitely generated submodules of $M$ are in
$\T_{\dim}$. The dimension of $M$ is equal to the supremum of the
dimensions of finitely generated submodules of $M$ by Cofinality.
But that supremum is 0, so $M$ is in $\T_{\dim}.$ \qed \end{pf}

\subsection{Comparing the Torsion Theories for $\ce$}
Let us summarize the situation with various torsion theories of a
ring $R$ from $\ce.$ The trivial torsion theory is the theory $(0,
\mathrm{Mod}_R),$ where Mod$_R$ is the class of all $R$-modules
and the improper torsion theory is the theory $(\mathrm{Mod}_R,
0).$ The various torsion theories for $R$ are ordered as follows:
\begin{center}
Trivial $\leq$ Classical = $\tau_{Q}$ = $(\smallt,\p)$ $\leq$
$(\T_{\dim}, \bigP_{\dim})$ = $(\T, \bigP)$ $\leq$ $(\bnd, \unb)$
$\leq$ Improper
\end{center}
where all of the above inequalities can be strict. The theory
$(\smallt,\p)$ can be nontrivial by Example 2.9 in \cite{Lu_book}.
The inequality $(\smallt,\p)\leq(\T, \bigP)$ can be strict by
Example 8.34 in \cite{Lu_book}. The inequality $(\T,
\bigP)\leq(\bnd, \unb)$ can be strict by example given in Exercise
6.5 in \cite{Lu_book}. Note that all of the above examples are
given for finite von Neumann algebras.

For any nontrivial $R$ the theory $(\bnd, \unb)$ is not improper
since $R$ is a module in $\unb.$

We have seen that the classes $\T$ and $\bnd$ coincide for
finitely generated modules. Next, we shall show that the classes
$\T$ and $\smallt$ coincide when restricted on the class of
finitely presented $R$-modules. First, we need the following
result proven in \cite{Lia1} for finite von Neumann algebras.

\begin{prop} If $P$ is finitely generated projective $R$-module,
then \[P\otimes_R Q=E(P).\] \label{EnvelopeIsTensor}
\end{prop}

For the proof, see Theorem 5.1 from \cite{Lia1}. In \cite{Lia1},
this result was proven for a finite von Neumann algebra but the
proof transfers to any $R\in \ce.$ The idea is to show that
$P\subseteq_e P\otimes_R Q$ (which holds because $Q$ is the
classical ring of quotients) and that $P\otimes_R Q$ is injective
(which holds because $P$ is a direct summand of some finitely
generated free $R$-module). The Lemma 5.1 from \cite{Lia1} also
holds for a ring in $\ce.$

Now we can prove the following.

\begin{prop} If $M$ is finitely presented $R$-module, then $\T M=\smallt M.$
\label{FinPresented}
\end{prop}
\begin{pf}
Since $\smallt M\subseteq \T M,$ it is sufficient to prove the
converse. We shall show that if $M$ has dimension zero (i.e. $M\in
\T$), then $M\otimes_R Q=0$ (i.e. $M\in \smallt$). Since $M$ is
finitely presented, there are finitely generated projective
modules $P_0$ and $P_1$ such that
\[0\rightarrow P_1\rightarrow P_0\rightarrow M\rightarrow0.\]

$\dim_R(M)=0$ and so $\dim_R(P_0)=\dim_R(P_1).$ Then
\[\dim_Q(P_0\otimes_R Q)=\dim_Q(E(P_0))=\dim_R(P_0)=\] \[=\dim_R(P_1)=\dim_Q(E(P_1))=\dim_Q(P_1\otimes_R Q)\]
by Proposition \ref{EnvelopeIsTensor} and part (1) of Lemma
\ref{ClosureAndSup}. $Q$ is $R$-flat, so
\[0\rightarrow P_1\otimes_R Q\rightarrow P_0\otimes_R Q\rightarrow M\otimes_R Q\rightarrow0\]
is an exact sequence. Thus, $\dim_Q(M\otimes_R
Q)=\dim_Q(P_0\otimes_R Q)-\dim_Q(P_1\otimes_R Q)=0.$  Moreover,
the modules $P_0\otimes_R Q$ and $P_1\otimes_R Q$ are finitely
generated projective $Q$-modules and hence $M\otimes_R Q$ is
finitely presented $Q$-module. But $Q$ is a regular ring so all
finitely presented modules are projective. Thus, $M\otimes_R Q$ is
a finitely generated projective module of dimension zero and so it
is trivial. \qed
\end{pf}

Before proving the main result of this subsection, let us prove
another corollary of Theorem \ref{DimensionInR}, Proposition
\ref{EnvelopeIsTensor} and Lemma \ref{ClosureAndSup}.

\begin{cor} For any $R$-module $M,$ \[\dim_Q(M\otimes_R Q)=\dim_R(M).\]
\label{dimQVSdimR}
\end{cor}
\begin{pf} If $M$ is a finitely generated projective module,
this follows from Lemma \ref{ClosureAndSup} and Proposition
\ref{EnvelopeIsTensor}. If $M$ is a submodule of any projective
$R$-module, write $M$ as a directed union its finitely generated
modules $M_i,$ $i\in I.$ All the modules $M_i$ are projective as
they are finitely generated submodules of a projective module and
$R$ is semihereditary. Thus,
\[\dim_R(M)=\sup_{i\in I} \dim_R(M_i)= \sup_{i\in I} \dim_Q(M_i\otimes_R
Q)=\dim_Q(M\otimes_R Q)\] by Cofinality of $\dim_R$ and $\dim_Q$
and the fact that tensor commutes with direct limit.

Finally, if $M$ is an arbitrary $R$-module, then $M$ is a quotient
of some projective module $P$ and its submodule $K$. Then,
$\dim_R(M)=\dim_R(P)-\dim_R(K)=\dim_Q(P\otimes_R
Q)-\dim_Q(K\otimes_R Q)=\dim_Q(M\otimes_R Q)$ by Additivity of
dimensions $\dim_R$ and $\dim_Q$ and since $Q$ is $R$-flat. \qed
\end{pf}

We now show how torsion theories reflect the ring-theoretic
properties of $R.$

\begin{thm}
\begin{enumerate}
\item $R$ is regular if and only if $(\smallt, \p)$ is trivial.
\item If $R$ is self-injective, then $(\T, \bigP)=(\bnd, \unb).$
\item The regular ring $Q$ of $R$ is semisimple if and only if
$(\smallt, \p)=(\T, \bigP)$ for $R.$

\item The following are equivalent
\begin{enumerate}
\item $R$ is semisimple,

\item $(\bnd, \unb)$ is trivial,

\item $(\T, \bigP)$ is trivial.
\end{enumerate}

\item $R$ is trivial if and only if $(\bnd, \unb)$ is improper.
\end{enumerate}
\label{TTEqual}
\end{thm}
\begin{pf} (1) $R$ is regular if and only if all the $R$-modules are
flat. But, $(\smallt, \p)$ is trivial if and only if all
$R$-modules are in $\p$ i.e. flat.

(2) If $R$ is self-injective, then $R=E(R).$ Thus, the torsion
theories cogenerated with $R$ and $E(R)$ coincide.

(3) Suppose that $Q$ is semisimple. We show that
$\T\subseteq\smallt$ by showing that every $R$-module $M$ with
dimension zero is such that $M\otimes_R Q=0.$ If $\dim_R(M)=0,$
then $\dim_Q(M\otimes_R Q)=\dim_R(M)=0$ by Corollary
\ref{dimQVSdimR}. Thus, $M\otimes_R Q$ is a projective (since $Q$
is semisimple) and of dimension 0. Hence, $M\otimes_R Q=0.$

Conversely, if $\smallt=\T,$ we shall show that every right ideal
$I$ of $Q$ is a direct summand (thus $Q$ is semisimple). Since
$\cl^Q_{\T}(I)$ is a direct summand by Proposition
\ref{VariousClosures}, it is sufficient to show that
$I=\cl^Q_{\T}(I).$ Let us look at $\cl^R_{\T}(I\cap R).$

$\cl^R_{\T}(I\cap R)/(I\cap R)=\T(R/(I\cap R))=\smallt(R/(I\cap
R))$ by assumption that $\smallt=\T.$ Thus, $\cl^R_{\T}(I\cap
R)\otimes_R Q =  (I\cap R)\otimes_R Q.$ Then,

\[\begin{array}{rcll}
\cl^Q_{\T}(I) & = & E(I) & (\mbox{by Proposition
\ref{VariousClosures}})\\
 & = & E(I\cap R)  & (\mbox{since }I\cap R\subseteq_e I)\\
 & = & E(\cl^R_{\T}(I\cap R)) & (\mbox{since }I\cap R\subseteq_e \cl^R_{\T}(I\cap R))\\
 & = & \cl^R_{\T}(I\cap R)\otimes_R Q & (\mbox{by Proposition
 \ref{EnvelopeIsTensor}})\\
 & =& (I\cap R)\otimes_R Q & (\mbox{by what we showed above})\\
 & \subseteq & I & (I\mbox{ is a right ideal}).
\end{array}\]

Since $I\subseteq \cl^Q_{\T}(I)$ always holds, $I=\cl^Q_{\T}(I).$

(4) (a) $\Rightarrow$ (b) If $R$ is semisimple, then all
nontrivial $R$-modules are projective and, thus, in $\unb.$ So,
$(\bnd, \unb)$ is trivial.

(b) $\Rightarrow$ (c) is trivial since $(\T, \bigP)\leq(\bnd,
\unb).$

(c) $\Rightarrow$ (a) If $(\T, \bigP)$ is trivial, then $(\smallt,
\p)$ is trivial and then $R$ is regular by (1). Thus, $R=Q.$ But
$(\T, \bigP)$ is trivial implies that $(\smallt, \p)=(\T, \bigP)$
and so $Q$ is semisimple by (3). Then, $R=Q$ is semisimple.

(5) If $(\bnd, \unb)$ is improper, then $R\cong\homo_R(R,R)=0.$
The converse is trivial. \qed
\end{pf}

Part (3) of the above theorem generalizes the result (Theorem 6.6)
from \cite{Tele_Teza} proven there for group von Neumann algebras.
Part (2) and Theorem \ref{dim=G=L} generalize Theorem 5.1 from
\cite{Tele_Teza}.

The order of torsion theories for $R$ implies that for every
$R$-module $M,$ we have a filtration:
\[\underbrace{0\subseteq\smallt}_{\smallt
M}\underbrace{M\subseteq\T }_{\T\p M} \underbrace{M\subseteq
M}_{\bigP M}\]

The quotient $\T M/\smallt M$ is equal to $\T\p M=\T(M/\smallt M)=
\cl_{\T}(\smallt M)/\smallt M$ since $\cl_{\T}(\smallt M) = \T M$
as one can easily show (see Proposition 4.3 in \cite{Lia1}).

For $M$ finitely generated, the finitely generated quotient $\p M$
splits as the direct sum of $\T\p M$ and $\bigP\p M = \bigP M$ and
thus we have a short exact sequence $ 0\rightarrow\smallt
M\rightarrow M\rightarrow \T\p M\oplus\bigP M \rightarrow 0.$ Of
course, $M$ also splits as $\T M\oplus\bigP M$ by Theorem
\ref{dim=G=L}.

If $Q$ is a regular ring from $\ce,$ then it is its own regular
ring (by Theorem \ref{RegularRing}) and is self-injective by
Proposition \ref{Q=max=classical}. By Theorem \ref{TTEqual}, the
torsion theories of $Q$ are ordered as follows:
\begin{center}
Trivial = $(\smallt,\p)$ $\leq$ $(\T, \bigP)$ = $(\bnd, \unb)$
$\leq$ Improper
\end{center}
The first inequality is strict if and only if $Q$ is not
semisimple and the second if and only if $Q$ is nontrivial. Thus,
for general regular ring in $\ce$, there seem to be only one
nontrivial and proper torsion theory of interest.

\section{Corollaries}
\label{SectionCorollaries}

\subsection{Finite Von Neumann Algebras.} If $A$ is a finite von
Neumann algebra, the center $Z$ of $A$ can be identified with
$C(X)$ (because $Z$ is the closed linear span of central
projections). Thus, the dimension $\dim_A$ is central-valued.

Also, a finite von Neumann algebra $A$ has a function $\tr_A:
A\rightarrow Z$ uniquely determined by the following properties:
\begin{itemize}
\item[(T1)] $\tr_A$ is $\Cset$-linear,

\item[(T2)] $\tr_A(ab)=\tr_A(ba),$

\item[(T3)] $\tr_A(c)=c$ for every $c\in Z,$

\item[(T4)] $\tr_A(a)$ is positive for every  positive $a$ (i.e.
$a=b^*b$ for some $b$).
\end{itemize}

The restriction of $\tr_A$ to the set of projections in $A$
satisfies the properties (D1) -- (D4) of Theorem
\ref{DimensionExists} (for proof and details see section 8.4 of
\cite{KadRingtwo}). Thus, we can use the center-valued trace to
define the center-valued dimension function $\dim_A$ using the
approach given in this paper.

In \cite{Lu_book}, \cite{Lu2} and \cite{Lia1} the real valued
dimension $\dim^{\Rset}_A$ of a module over a finite von Neumann
algebra $A$ was considered and results analogous to those we
proved here for a ring in $\ce,$ were proven for $A$. Since the
real-valued dimension $\dim^{\Rset}_A$ depends on the
complex-valued trace used (that is not unique), it was surprising
that the torsion theory $(\T_{\dim^{\Rset}_A},
\bigP_{\dim^{\Rset}_A})$ coincided with Lambek and Goldie theories
(Proposition 4.2 in \cite{Lia1}) regardless of the complex-valued
trace used to define the dimension.

A direct corollary of our Theorem \ref{dim=G=L} and Proposition
4.2 in \cite{Lia1} is the following
\begin{cor}
\[(\T_{\dim_A}, \bigP_{\dim_A})=(\T_{\dim^{\Rset}_A}, \bigP_{\dim^{\Rset}_A}).\]

i.e. for every $A$-module $M$, the central-valued dimension of $M$
is zero if and only if a real-valued dimension of $M$ is zero.
\label{CorollaryDimensions}
\end{cor}

It is interesting to note that in the case of a finite von Neumann
algebra $A$, the algebra of affiliated operators $U$ does not
automatically come equipped with a trace (and thus a dimension)
function since $U$ might not be a finite von Neumann algebra. In
\cite{Lu_book} it is shown that we can still get the real-valued
$U$-dimension of any $U$-module in the following way: if $S$ is a
finitely generated projective $U$-module, the $U$-dimension of $S$
is the $A$-dimension of a finitely generated projective $A$-module
$P$ such that $P\otimes_A U\cong S$ (this is well defined by
Theorem 8.22, \cite{Lu_book}). Then one proves that $U$ satisfies
Theorem \ref{LueckDimension}. So, the $U$-dimension can be
extended to all $U$-modules (Lemma 8.27, \cite{Lu_book}).
Moreover, Corollary \ref{dimQVSdimR} holds for $A$ and $U$
(Theorem 8.29, \cite{Lu_book}).

The regular ring $Q$ of a ring $R$ from the class $\ce$ is
automatically equipped with a dimension since $Q$ is also in $\ce$
(see subsection \ref{Matrices}). Thus, by the results of this
paper, it readily follows that $U$ has the dimension function
$\dim_U$ satisfying Theorem \ref{LueckDimension} and Corollary
\ref{dimQVSdimR}.

Moreover, $Q$ has all the properties of $R$ and more (regularity,
self-injectiveness etc). In contrast, the algebra of affiliated
operators $U$ of a finite von Neumann algebra $A$ is regular but
might not be a von Neumann algebra (i.e. $U$ may loose some
properties of $A$).

\subsection{Cofinal-measurable Modules.}
Using the dimension function, we can view the theory $(\smallt,
\p)$ for $R\in\ce$ from a different angle. We say that an
$R$-module $M$ is {\em measurable} if it is a quotient of a
finitely presented module of dimension zero. $M$ is {\em
cofinal-measurable} if each finitely generated submodule is
measurable. The class of cofinal-measurable modules is equal to
the class $\smallt$ by Lemma 8.36 (2) from \cite{Lu_book}. The
proof there is given for a group von Neumann algebra but it
converts to any $R$ from $\ce.$

\subsection{$K_0$-theorem.}

Theorem 9.20 (1) from \cite{Lu_book} states that the inclusion $i:
A\rightarrow U$ of a finite von Neumann algebra $A$ to its algebra
of affiliated operators $U$ (same as regular ring of $A$) induces
an isomorphism of monoids $\mu :\mathrm{Proj}(A)$ $\rightarrow$
$\mathrm{Proj}(U)$ given by $[P]\mapsto[P\otimes_{A}U]$ and an
isomorphism $\mu: K_0(A)\rightarrow K_0(U).$

The proof of this theorem relies solely on Theorem 8.22 in
\cite{Lu_book}. Theorem 8.22 holds for a ring $R\in \ce$ and its
regular ring $Q.$ Thus, the result on $K_0$-theories holds for
$R\in \ce$ and its regular ring $Q.$

In \cite{Lia1}, the inverse of the isomorphism $\mu$ is described.
Namely, the following corollary is proven for a finite von Neumann
algebras. The proof transfers, word-for-word, to a ring from class
$\ce.$

\begin{cor}
For every finitely generated projective $R$-module $M$, there is
an one-to-one correspondence \[\{\mbox{direct summands of }M\}
\longleftrightarrow\{\mbox{direct summands of }E(M)\}\] given by
$K \mapsto$ $K\otimes_{R}Q = E(K).$ The inverse map is given by
$L\mapsto L\cap M.$  This correspondence induces an isomorphism of
monoids $\mu :\mathrm{Proj}(R)$ $\rightarrow$ $\mathrm{Proj}(Q)$
and an isomorphism \[\mu: K_0(R)\rightarrow K_0(Q)\] given by
$[P]\mapsto[P\otimes_{R}Q]$ with the inverse $[S]\mapsto[S\cap
R^n]$ if $S$ is a direct summand of $Q^n.$ \label{moja K 0
theorem}
\end{cor}

Since we have Corollary \ref{Johnson}, Proposition
\ref{EnvelopeIsTensor}, and Theorem \ref{SplittingOfBnd} (recall
that complements are $(\T,\bigP)$-closed modules and
$(\T,\bigP)=(\bnd, \unb)$ for finitely generated modules), it is
not hard to see why Corollary \ref{moja K 0 theorem} holds for
$R\in \ce.$

Let us mention that Handelman proved (Lemma 3.1 in \cite{Handel})
that for every finite Rickart $C^*$-algebra $A$ such that every
matrix algebra over $A$ is also Rickart, the inclusion of $A$ into
its regular ring $Q$ induces an isomorphism $\mu:
K_0(A)\rightarrow K_0(Q).$ By Theorem 3.4 in \cite{AG}, a matrix
algebra over a Rickart $C^*$-algebra is a Rickart $C^*$-algebra.
Thus, $K_0(A)$ is isomorphic to $K_0(Q)$ for every finite Rickart
$C^*$-algebra. Corollary \ref{moja K 0 theorem} describes the
inverse of this isomorphism.

\subsection{Questions.} We conclude by listing some interesting
questions.
\begin{enumerate}
\item Is it possible to obtain the same results using the weaker
axioms than (A1) -- (A9)? Note that (A8) and (A9) are particularly
strong.

\item Let $R$ be in $\ce$ and $Q$ be the regular ring of $R$. Does
$Q$ semisimple imply $R$ semisimple? Note that this is the case
for group von Neumann algebras by Exercise 9.11 from
\cite{Lu_book} (if the algebra of affiliated operators of a group
von Neumann algebra $\vng$ is semisimple, then the group $G$ is
finite so $\vng$ is finite dimensional over $\Cset$ and, thus,
semisimple).

\item Let $R$ be in $\ce$ and $(\T, \bigP)=(\bnd, \unb)$ for $R$.
Is $R$ self-injective? The converse holds by Theorem
\ref{TTEqual}.
\end{enumerate}

\end{document}